\documentclass{amsart}

\usepackage[utf8]{inputenc}
\usepackage{amsmath,amssymb,amsthm}
\usepackage{dsfont} 
\usepackage[a4paper, total={5.5in, 8.5in}]{geometry}

\usepackage{hyperref}
\usepackage{cleveref}
\usepackage[style=alphabetic,sorting=anyt,sortcites=true,backend=biber,maxbibnames=5,maxcitenames=5]{biblatex}
\newtheorem{thm}{Theorem}[section]
\newtheorem{lem}[thm]{Lemma}
\newtheorem{cor}[thm]{Corollary}
\newtheorem{prop}[thm]{Proposition}
\newtheorem{ex}[thm]{Example}

\newtheorem{con}[thm]{Conjecture}
\newtheorem*{prob*}{Open problem}

\theoremstyle{definition}

\newtheorem{defi}[thm]{Definition}

\theoremstyle{remark}

\newtheorem{rem}[thm]{Remark}
\newtheorem*{rem*}{Remark}

\DeclareMathOperator{\End}{End}
\DeclareMathOperator{\Aut}{Aut}

\DeclareMathOperator{\tr}{tr}
\DeclareMathOperator{\GL}{GL}
\DeclareMathOperator{\Spec}{Spec}
\DeclareMathOperator{\ESpec}{ESpec}
\DeclareMathOperator{\Fix}{Fix}

\newcommand{\NN}{{\mathbb{N}}}
\newcommand{\ZZ}{{\mathbb{Z}}}
\newcommand{\QQ}{{\mathbb{Q}}}
\newcommand{\RR}{{\mathbb{R}}}
\newcommand{\CC}{{\mathbb{C}}}

\newcommand{\lie}[1]{{\mathfrak {#1}}}

\begin{document}
\title[Fixed points of diffeomorphisms of nilmanifolds]{Fixed points of diffeomorphisms on nilmanifolds with a free nilpotent fundamental group} 

\author[K. Dekimpe]{Karel Dekimpe}\thanks{Research supported by long term structural funding -- Methusalem grant of the Flemish Government.}
\author[S. Tertooy]{Sam Tertooy}
\address{KU Leuven Campus Kulak Kortrijk\\
E. Sabbelaan 53\\
8500 Kortrijk\\
Belgium}
\email{Karel.Dekimpe@kuleuven.be}
\email{Sam.Tertooy@kuleuven.be}

\author[A. R. Vargas]{Antonio R. Vargas}
\address{KU Leuven\\
Department of Mathematics\\
Celestijnenlaan 200B\\
3000 Leuven\\
Belgium}
\email{antonio.vargas@kuleuven.be}

\date{\today}

\maketitle

	\begin{center}
		This is an Accepted Manuscript of an article published by International Press of Boston in The Asian Journal of Mathematics on 21 Aug 2020, available online:  \href{https://doi.org/10.4310/AJM.2020.v24.n1.a6}{https://doi.org/10.4310/AJM.2020.v24.n1.a6}.
\end{center}

\begin{abstract}
	Let $M$ be a nilmanifold with a fundamental group which is free $2$-step nilpotent on at least 4 generators. We will show that for any nonnegative integer $n$ there exists a self-diffeomorphism $h_n$ of $M$ such that $h_n$ has exactly $n$ fixed points and any self-map $f$ of $M$ which is homotopic to $h_n$ has at least $n$ fixed points.  We will also shed some light on the situation for less generators and also for higher nilpotency classes.
\end{abstract}

\section{Reidemeister-Nielsen fixed point theory}\label{intro}
The aim of this paper is to study fixed points of self-homeomorphisms of nilmanifolds. In this first section, we will recall the basics of Reidemeister-Nielsen theory.

\medskip

Let $f:X\to X$ be a selfmap of a closed manifold $X$ (or even a compact polyhedron). Let $\Fix(f)=\{x\in X\;|\; f(x) =x\}$ denote the set of fixed points of $f$. The aim of Reidemeister-Nielsen theory is to obtain a good estimate for the minimal size of $\Fix(g)$ for all $g$ homotopic to a given $f$. A first result in this direction is given by the Lefschetz number $L(f)$ of $f$, this number is given by 
\[ L(f) =\sum_{i=0}^{\dim X} (-1)^i \tr (f_{\ast,i}: H_i(X,\QQ) \to H_i(X,\QQ)\,).\]
Lefschetz' theorem says that if $L(f)\neq 0$, then $f$ must have a fixed point. As homotopic maps induce the same morphisms on homology, we can also deduce that if $L(f)\neq 0$, then any map $g$ homotopic to $f$ has a least one fixed point. Unfortunately, Lefschetz' theorem does not give any information on the exact amount of fixed points and says nothing in case $L(f)=0$. One can associate a more sophisticated number to $f$, namely the Nielsen number $N(f)$, which in many cases gives full information.

\medskip

Let us spend some time in explaining how this Nielsen number is defined. For details we refer the reader to \cite{jm06-1,jian83-1}. Consider the universal covering space $p: \tilde{X} \to X$ of $X$. For any lift $\tilde{f}:\tilde{X} \to \tilde{X}$ of $f$ we have that 
$p(\Fix(\tilde{f}))\subseteq \Fix(f)$. In fact 
\[ \Fix(f) = \bigcup_{\tilde{f}} p(\Fix(\tilde{f})),\]
where the union is taken over all possible lifts $\tilde{f}$ of $f$. We call two lifts $\tilde{f}_1$ and $\tilde{f}_2$ equivalent (notation $\tilde{f}_1 \sim \tilde {f}_2$) if $\tilde{f}_1 = \gamma \circ \tilde{f}_2 \circ \gamma^{-1}$ for some covering transformation $\gamma$.
It then follows that 
\begin{itemize}
\item $\tilde{f}_1 \sim \tilde{f}_2 \Rightarrow p\left( {\rm Fix}(\tilde{f}_1)\right)=p\left( {\rm Fix}(\tilde{f}_2)\right)$,
\item $\tilde{f}_1 \not \sim \tilde {f}_2 \Rightarrow p\left( {\rm Fix}(\tilde{f}_1)\right) \cap p\left( {\rm Fix}(\tilde{f}_2)\right)=\emptyset$.
\end{itemize}
For any equivalence class $[\tilde{f}]$ of lifts we call $p(\Fix(\tilde{f}))$ the fixed point class determined by $[\tilde{f}]$. The above shows that this fixed point class is independent of the chosen representative and that different equivalence classes determine disjoint fixed point classes. Hence $\Fix(f)$ is the disjoint union of its fixed point classes. Note that it is possible that some fixed point classes are empty (but still they are considered to be a fixed point class!). The Reidemeister number of the map $f$, denoted by $R(f)$ is the number of fixed point classes of $f$ (or the number of equivalence classes of lifts of $f$). 

\medskip

There is an algebraic way to determine the Reidemeister number of a map $f$, using the notion of twisted conjugacy. Let $\Gamma$ be any group and suppose that $\varphi:\Gamma \to \Gamma$ is an endomorphism. We will say that two elements $\gamma$ and $\gamma'$ of $\Gamma$ are twisted conjugate (or Reidemeister equivalent) if and only if there exists a
$\mu \in \Gamma$ such that $\gamma = \mu \gamma \varphi(\mu)^{-1}$. Being twisted conjugate is an equivalence relation on $\Gamma$ and the number of equivalence classes is denoted by $R(\varphi)$ (and is called the Reidemeister number of $\varphi$).

\medskip

To explain the connection with the Reidemeister number of a map $f$, we fix a lift (to the universal covering space) $\tilde{f}_0$  of $f$ and let $\Gamma$ be the group of covering transformations of the universal covering $p:\tilde{X} \to X$. Recall that $\Gamma$ is isomorphic to the fundamental group of $X$. Any other lift $\tilde{f}$ of $f$ can be written uniquely as a composition 
$\tilde{f}= \gamma\circ \tilde{f}_0$ and conversely any map $\gamma \circ \tilde{f}_0$ is a lift of $f$. Hence there is a 
one-to-one correspondence between lifts of $f$ and the elements of $\Gamma$. There is a morphism 
\[ f_\ast: \Gamma\to \Gamma: \gamma \mapsto f_\ast(\gamma)  \mbox{ determined by } f_\ast(\gamma) \circ \tilde{f}_0 = 
\tilde{f}_0\circ \gamma.\]
This morphism depends on the choice of $\tilde{f}_0$, but only up to an inner automorphism of $\Gamma$, and in fact, for the good choices of base points it coincides with the morphism induced by $f$ on the fundamental group of $X$. 

It is well known (and easy to check) that 
\[ \gamma \circ \tilde{f}_0 \sim \gamma' \circ \tilde{f}_0 \Leftrightarrow \exists \mu \in \Gamma: \gamma = 
\mu \circ \gamma' \circ f_\ast (\mu)^{-1}.\]
It follows that equivalence of  two lifts determined by elements $\gamma$ and $\gamma'$ corresponds exactly to the fact that these two elements are twisted conjugate (using the morphism $f_\ast$) and hence fixed point classes are in one-to-one correspondence with twisted conjugacy classes of $f_\ast$, i.e.\ $R(f)=R(f_\ast)$.

\medskip

One of the main ingredients of Nielsen theory is the so-called fixed point index. There is an axiomatic way to attach to each fixed point class ${\mathbf F}$ of a map $f$ an integer $I(f,{\mathbf F})$. It is not easy to explain in a few words how this index is determined, but the idea is that if the index is non-zero, then this fixed point class cannot disappear (become empty)  under a homotopy. Therefore, a fixed point class ${\mathbf F}$ with $I(f,{\mathbf F})\neq 0$ is called an essential fixed point class. To have at least some idea about the meaning of the index, we consider the case when $X$ is a differentiable manifold and ${\mathbf F}=\{x_0\}$ is a fixed point class consisting of an isolated fixed point $x_0$. In this case 
\[ I(f,{\mathbf F})={ \rm sgn} \det(\mathds{1} - df_{x_0})\]
where $df_{x_0}:T_{x_0}X\to T_{x_0}X$ is the differential of $f$ at $x_0$, $\mathds{1}$ denotes the identity map of $T_{x_0}X$ and 
for $r\in \RR$ we let sgn$(r)=-1, 0$ or $1$ when $r<0$, $r=0$ and $r>0$ respectively. The Nielsen number of the map $f$ is then  defined as the number of essential fixed point classes:
\[ N(f) = \# \{ {\mathbf F}\;|\; {\mathbf F}\mbox{ is a fixed point class with } I(f,{\mathbf F})\neq 0\}.\]
Obviously, $N(f)$ is a lower bound for the number of fixed points of $f$ and as $N(f)$ is a homotopy invariant, we also have that 
\[ N(f)\leq \# \Fix(g)\mbox{ for all } g\sim f.\]
In fact, Wecken showed that in many cases $N(f)$ is a sharp lower bound for the minimal number of fixed points of all maps $g$ which are homotopic to $f$.
\begin{thm}[Wecken \cite{weck42-1}]
Let $X$ be a compact manifold, possibly with boundary, of dimension $\geq 3$ and $f:X\to X$ be a self map of $X$. Then there exists a map $g:X\to X$ homotopic to $f$ such that $N(f)=\# \Fix(g)$.
\end{thm}

\section{Nilmanifolds and their Reidemeister spectrum}

From now on we will focus on nilmanifolds. 
\begin{defi} Let $G$ be a connected and simply connected nilpotent Lie group.
A lattice of $G$ is a discrete and cocompact subgroup $N$ of $G$. Here cocompact means that the quotient space $N\backslash G$, where $N$ is acting by left translations on $G$, is compact.\\ 
A nilmanifold is a quotient manifold $N\backslash G$, where $N$ is a lattice of some connected and simply connected nilpotent Lie group $G$.
\end{defi}

\medskip

The following results due to Mal'cev are well known (\cite{malc51-1, ragh72-1})
\begin{itemize}
\item If $N$ is a lattice of a connected and simply connected nilpotent Lie group, then $N$ is a finitely generated torsion-free nilpotent group.
\item Conversely, if $N$  is a finitely generated torsion-free nilpotent group, then there exists a connected and simply connected nilpotent Lie group $G$ and an embedding $i:N\to G$, such that $i(N)$ is a lattice of $G$. Moreover, $G$ is unique up to isomorphism. We refer to $G$ as being the Mal'cev completion of $N$. In practice, we do not write the $i$ and we assume that we know how $N$ sits inside $G$ as a lattice.
\item If $\varphi\in \End(N)$ is an endomorphism of a finitely generated torsion-free nilpotent group $N$, then $\varphi$ extends uniquely  to a Lie group endomorphism $\tilde{\varphi}:G\to G$ of the Mal'cev completion $G$ of $N$. Moreover, if $\varphi$ is an automorphism, then $\tilde{\varphi}$ is also an automorphism.
\end{itemize}

\medskip

The Reidemeister-Nielsen fixed point theory is very well developed for nilmanifolds. We wil explain some details below.

\medskip

We fix a connected and simply connected nilpotent Lie group $G$ and denote its Lie algebra by $\lie{g}$. Recall that for simply connected nilpotent Lie groups, the exponential map $\exp:\lie{g} \to G$ is a diffeomorphism and we will use $\log$ to denote its inverse. For any Lie group endomorphism $\varphi$ of $G$, its differential $d\varphi$ induces a Lie algebra endomorphism of $\lie{g}$ and conversely,  any endomorphism of $\lie{g}$ can be seen as being induced by a Lie group endomorphism. Moreover, $\varphi$ is an automorphism if and only if $d \varphi$ is.   Now consider a lattice $N$ of $G$ and a self-map $f$ of the nilmanifold $N\backslash G$. Note that $G$ is the universal covering space of $N\backslash G$ and $N$ is (isomorphic to) the group of covering transformations of the covering $p:G \to N\backslash G$. Assume that $f$ induces the endomorphism $f_\ast=\varphi:N\to N$ (as in section~\ref{intro}). By abuse of notation, we also use $\varphi$ to denote the unique extension of $\varphi$ to $G$ and so we have an induced endomorphism $d\varphi$ on $\lie{g}$. We then have (\cite{anos85-1,jm06-1}):
\begin{equation}\label{anosov}
N(f)=|\det(\mathds{1}-d\varphi)|.
\end{equation}
Moreover, on nilmanifolds there is a strong relation between the Reidemeister number of a map and the Nielsen number 
(\cite{hk97-1}): If $N(f)\neq 0$, then $N(f)=R(f)$ and $N(f)=0\Leftrightarrow R(f)=\infty$. In fact, more can be said: on a nilmanifold it holds that all fixed point classes have the same index, which is either $\pm 1$ or $0$ (see also \cite{fl13-1, kll05-1}). So either all fixed point classes are essential or all of them are inessential.

\medskip

There is an algebraic way to construct maps on a given nilmanifold $N\backslash G$. Choose any endomorphism $\varphi$ of $G$ with $\varphi(N)\subseteq N$. (This means that $\varphi$ is the unique extension of an endomorphism on $N$). Then $f_\varphi:N\backslash G \to N\backslash G:\; N x \mapsto N \varphi(x)$ is a well defined map and moreover $f_{\varphi\ast} =\varphi$ by construction. The map $f_\varphi$ is a diffeomorphism if and only if $\varphi(N)=N$.  We will refer to $f_\varphi$ as a nilmanifold endomorphism and if $\varphi(N)=N$ as a nilmanifold automorphism. Conversely, if $f$ is any self-map of the nilmanifold $N\backslash G$, we can consider $\varphi=f_\ast$ and if we let $g$ denote the nilmanifold endomorphism determined by $\varphi$ we have that $f_\ast=g_\ast=\varphi$. From this, it follows that $f\sim g$. We can conclude that any self-map of a nilmanifold is homotopic to a nilmanifold endomorphism.

\medskip

The following result can be found (even in a more general setting)  in \cite[Remark 9.4]{fl13-1}:
\begin{prop}\label{weckenmaps} Let $M=N\backslash G$ be a nilmanifold and let $\varphi$ be an endomorphism of $N$ determining a nilmanifold endomorphism $f_\varphi:M\to M$. 
If $N(f_\varphi)\neq 0$, then $\# \Fix(f_\varphi)=N(f_\varphi)$. 
\end{prop}
The above proposition says that when $N(f_\varphi)\neq 0$ then any map homotopic to $f_\varphi$ has at least as many fixed points as the map $f_\varphi$. 

\medskip

Let $\Gamma$ be any group. The Reidemeister spectrum of $\Gamma$ is the set 
\[ \Spec_R(\Gamma)=\{ R(\varphi) \; |\; \varphi\in \Aut(\Gamma) \}.\]
So we have that $\Spec_R(\Gamma)\subseteq \NN \cup \{\infty\}$ (where we use $\NN$ to denote the set of positive integers). 
When $\Spec_R(\Gamma)= \NN \cup \{\infty\}$ we say that $\Gamma$ has full Reidemeister spectrum. A group with $\Spec_R(\Gamma)=\{\infty\}$ is said to have the $R_\infty$ property. 

From the facts we described above we more or less immediately get the following result.
\begin{thm}\label{thm:main}
If $N\backslash G$ is a nilmanifold and $n$ is a positive integer such that $n\in \Spec_R(N)$, then $N\backslash G$ admits a self-diffeomorphism $h_n$ with exactly $n$ fixed points and moreover any self-map $f$ of $N\backslash G$ which is homotopic to $h_n$ has at least $n$ fixed points.
\end{thm}
\begin{proof}
As $n\in \Spec_R(N)$, there exists an automorphism $\varphi\in \Aut(N)$ with $R(\varphi)=n$. If we now take for $h_n$ the self-diffeomorphism $f_\varphi$ then $R(h_n)=R(f_\varphi)=R(\varphi)=n$ and since $n\neq \infty$ we also have that $N(h_n)=R(h_n)=n$. Therefore \cref{weckenmaps} says exactly that $h_n$ has $n$ fixed points and moreover  any map homotopic to $h_n$ has at least $N(h_n)=n$ fixed points.
\end{proof}

\begin{cor}\label{cor:main} If $N\backslash G$ is a nilmanifold and $\Spec_R(N)$ is full, then for all non-negative integers $n$ there exists
 a self-diffeomorphism $h_n$ of $N\backslash G$ with exactly $n$ fixed points and moreover any self-map $f$ of $N\backslash G$ which is homotopic to $h_n$ has at least $n$ fixed points.
\end{cor}
\begin{proof}
The only thing we still have to show is the existence of a self-diffeomorphism $h_0$ of $N\backslash G$ without fixed points. This is easy: choose any element $z\in G$ such that $z\in Z(G)$ (the centre of $G$) and $z\not\in N$. This is possible because $N$ is a discrete subgroup of  $G$. Now, the map $h_0: N\backslash G \to N \backslash G: N x \mapsto Nzx $ is a well defined diffeomorphism without fixed points. 
\end{proof}

\section{Free nilpotent groups}

For any group $\Gamma$ we define the terms $\gamma_i(\Gamma)$ of the lower central series of $\Gamma$ inductively by 
$\gamma_1(\Gamma)=\Gamma$ and $\gamma_{i+1}(\Gamma)=[\Gamma,\gamma_i(\Gamma)]$. Analogously, for a Lie algebra $\lie{g}$ we have that $\gamma_1(\lie{g}) = \lie{g}$ and $\gamma_{i+1}(\lie{g}) = [\lie{g},\gamma_i(\lie{g})]$.

From now onwards, we concentrate on the case where $N$ is a free nilpotent group. Fix an integer $r\geq 2$ and denote by $F_r$ the free group on $r$ generators. The free nilpotent group of class $c$ on $r$ generators is the group 
\[ N_{r,c}= \frac{F_r}{\gamma_{c+1} (F_r)} .\]
Let $G_{r,c}$ be the Mal'cev completion of $N_{r,c}$ and $\lie{g}_{r,c}$ be the Lie algebra corresponding to $G_{r,c}$. Then 
$\lie{g}_{r,c}$ is the free $c$-step nilpotent Lie algebra on $r$ generators. 

For the free nilpotent group $N_{r,c}$, it is known that the terms of the lower central series $\gamma_i(N_{r,c})$ coincide with those of the upper central series (e.g.\ $\gamma_c(N_{r,c})=Z(N_{r,c})$, the centre of $N_{r,c}$). It follows that 
all quotients $\gamma_{i}(N_{r,c})/\gamma_{i+1}(N_{r,c})$ are free abelian groups.

\medskip

Let $f$ be a map on the nilmanifold $N_{r,c}\backslash G_{r,c}$ inducing a morphism $\varphi\in \End(N_{r,c})$. As the terms of the lower central series of a group are fully characteristic, $\varphi$ induces morphisms
\[ \varphi_i: \gamma_{i}(N_{r,c})/\gamma_{i+1}(N_{r,c}) \to \gamma_{i}(N_{r,c})/\gamma_{i+1}(N_{r,c}): n \gamma_{i+1}(N_{r,c})\mapsto \varphi(n) \gamma_{i+1}(N_{r,c}).\]
Let $a_1, a_2, \ldots, a_{k_i}\in  \gamma_{i}(N_{r,c})$ so that their canonical projections 
$\bar{a}_1, \bar{a}_2, \ldots, \bar{a}_{k_i}$ (where $\bar{a}_j = a_j\gamma_{i+1}(N_{r,c})$) form a free generating set of the free abelian group $\gamma_{i}(N_{r,c})/\gamma_{i+1}(N_{r,c})\cong\ZZ^{k_i}$. With respect to these generators the map $\varphi_i$ can be expressed by a matrix $M_i\in \ZZ^{k_i\times k_i}$.

\medskip

We can do the same on the Lie algebra level. The map $d\varphi$ induces morphisms (linear maps) 
\[ d \varphi_i :\gamma_i(\lie{g}_{r,c})/\gamma_{i+1}(\lie{g}_{r,c}) \to \gamma_i(\lie{g}_{r,c})/\gamma_{i+1}(\lie{g}_{r,c}): 
X+  \gamma_{i+1}(\lie{g}_{r,c})\mapsto d\varphi(X)+ \gamma_{i+1}(\lie{g}_{r,c}). \]
We can now rewrite \eqref{anosov} to obtain
\begin{equation}\label{anosovprod}
N(f)=|\det(\mathds{1}-d\varphi)|=\prod_{i=1}^c |\det(\mathds{1}-d\varphi_i)|.
\end{equation}
Let $a_1, \ldots , a_{k_i}$ be as above and take $A_j=\log(a_j)$. Then the natural projections $\bar{A}_1$, $\bar{A}_2$, $\ldots$,
$\bar{A}_{k_i}$ form a basis of $\gamma_i(\lie{g}_{r,c})/\gamma_{i+1}(\lie{g}_{r,c})$ and when we express $d\varphi_i$ with respect to this basis we find exactly the same matrix $M_i$. So we can also write \eqref{anosovprod} as
\begin{equation}\label{anosovprod2}
N(f)=|\det(\mathds{1}-d\varphi)|=\prod_{i=1}^c |\det(\mathds{1}-d\varphi_i)|=\prod_{i=1}^c |\det(\mathds{1}-\varphi_i)|.
\end{equation}

\medskip

In what follows we will need a so-called Hall basis of a free (nilpotent) Lie algebra. Let us recall the construction of such a basis. More details (and a more formal treatment) can be found in e.g.\ \cite[Chapter IV]{serr92-1}. Assume that $\lie{f}_r$ is a free Lie algebra generated by $r$ generators $X_1,\; X_2, \ldots, X_r$. A hall basis $H$ of $\lie{f}_r$ is a totally ordered set and a vector space basis of $\lie{f}_r$ which is built up recursively as a union $\displaystyle H=\bigcup_{n\in \NN} H_n$ where $H_n$ consists of basis vectors of length $n$ (these are $n$-fold Lie brackets in the generators) according to the following rules:

\begin{itemize}
\item $H_1=\{X_1, X_2, \ldots, X_r\}$,  and we order these basis vectors of length $1$ as follows
$X_1 < X_2 < \cdots < X_r$. 
\item Assume that $n\geq 2$ and $H_k$ has been defined for all $k<n$ and the order is already given 
on the set $\displaystyle \bigcup_{1\leq k \leq n-1} H_k$, then the elements $X$ of $H_n$ of length $n$ are formed as a Lie bracket of the form $X=[U,V]$ where $U\in H_k$ and $V\in H_l$ with $k+l=n$ and such that the following conditions are satisfied:
\begin{enumerate}
\item $U<V$, and
\item if $V=[V_1,V_2]$ for some $V_1\in H_{l_1}$ and $V_2\in H_{l_2}$, then $V_1\leq U$.
\end{enumerate}
\item The order on $\displaystyle \bigcup_{1\leq k \leq n-1} H_k$ is extended to an order on $\displaystyle \bigcup_{1\leq k \leq n} H_k$, by requiring that elements of length $\leq n-1$ are smaller than elements of length $n$ and by choosing any total order on the elements of $H_n$.
\end{itemize}

\begin{rem} The canonical images of the elements of length $\leq c$ (so of $H_1\cup H_2 \cup \cdots \cup H_c)$ form a basis of $\lie{g}_{r,c}$.
\end{rem}

\begin{ex}\label{exhall} The elements of $H_2$ are the elements of the form $[X_i,X_j]$ with $1 \leq i < j \leq r$.\\
The elements of $H_3$ are those of the form 
\[ [X_i,[X_j,X_k]] \mbox{ where $1 \leq j < k \leq n $ and $1 \leq j \leq i \leq r$}.\]  
The elements of $H_4$ depend on the choice of ordering we took for the elements of length 2.
\end{ex}

According to \eqref{anosovprod2}, to study fixed points of diffeomorphisms (or homeomorphisms or any map) of the nilmanifold $N_{r,c}\backslash G_{r,c}$, it is essential for us to be able to compute the determinants $\det( \mathds{1}-\varphi_i)$ for a given morphism $\varphi$ of $N_{r,c}$. This is equivalent to understanding the eigenvalues of the morphisms $\varphi_i$. The lemma below (which is a more explicit version of \cite[Lemma 2.4]{dg14-1}) shows that these are completely determined by the eigenvalues of $\varphi_1$. In order to be able to formulate this lemma, we need to introduce a map from a Hall basis to the complex numbers.

\begin{defi}
Fix an $r$-tuple of complex numbers $\lambda=(\lambda_1,\lambda_2, \ldots, \lambda_r)$. Let $H$ be a Hall basis of the free Lie algebra $\lie{f}_r$. We define a map 
$ f_\lambda : H \to \CC $  
recursively by 
\begin{itemize}
\item $\forall i \in \{1,2, \ldots , r\}$: $f_\lambda(X_i) = \lambda_i$.
\item Let $n\geq 2$ and assume that $f_\lambda(X)$ has been defined for all $X\in H_k$ with $1\leq k \leq n-1$. Now consider $X\in H_n$, then $X=[U,V]$ for some $U\in H_k$ and $V\in H_l$ with $k+l=n$. We let $f_\lambda(X)=f_\lambda(U) f_\lambda(V)$.
\end{itemize}
We will say that $f$ is determined by $\lambda$.
\end{defi}
\begin{ex} \label{exlambda}
We have that 
\[ f_\lambda([X_i,X_j]) = \lambda_i \lambda_j \mbox{ and } f_\lambda([X_i,[X_j,X_k]]=\lambda_i \lambda_j \lambda_k.\]
\end{ex}

\begin{lem}\label{lem:eigenvalues}
Let $r\geq 2$ and $c\geq1$ be positive integers and assume that $\varphi \in \End(N_{r,c})$ is a morphism inducing  morphisms $\varphi_i$ on the quotients $\gamma_{i}(N_{r,c})/\gamma_{i+1}(N_{r,c})$ ($1\leq i \leq c$). \\
Let $\lambda_1,\;\lambda_2,\; \cdots, \lambda_r$ be the eigenvalues of $\varphi_1$ (each eigenvalue is listed as many times as its multiplicity). Let $H$ be a Hall basis of the free Lie algebra $\lie{f}_r$ and $\lambda=(\lambda_1,\lambda_2,\ldots, \lambda_r)$. Let $f_\lambda:H\to \CC$ be the map associated to $\lambda$. Then the eigenvalues of $\varphi_i$ ($1\leq i \leq c$) are given as the collection of values 
\[ f_\lambda (X) \mbox{ where $X$ ranges over all elements of $H_i$}.\]
In this way each eigenvalue is then listed as many times as its multiplicity.
\end{lem}

\begin{proof}
Let $d\varphi$ denote the corresponding morphism on the Lie algebra $\lie{g}_{r,c}$. As mentioned before,
the eigenvalues of $\varphi_i$ are the same as the eigenvalues of $d\varphi_i$, the morphism induced by $d \varphi$ on $\gamma_i(\lie{g}_{r,c})/\gamma_{i+1}(\lie{g}_{r,c})$ (as they can be represented by the same matrix). 
It is well known that the semisimple part of $d\varphi$ is also an automorphism of  $\lie{g}_{r,c}$ (See e.g.\ \cite[Corollary 2, page 135]{sega83-1}) having the same eigenvalues as $d\varphi$ (also on each quotient $\gamma_i(\lie{g}_{r,c})/\gamma_{i+1}(\lie{g}_{r,c})$). Therefore, we may assume that $d\varphi$ is semisimple. Let $\lie{g}^\CC_{r,c}=\lie{g}_{r,c}\otimes_\RR \CC$  be the complexification of $\lie{g}_{r,c}$, then there exists a basis of  $\lie{g}_{r,c}^\CC$ consisting of eigenvectors for $d\varphi$ (which we can also consider as being a morphism of $\lie{g}_{r,c}^\CC$). It follows that we can find $r$ eigenvectors $X_1,X_2, \ldots, X_r$ of $\lie{g}_{r,c}^\CC$ such that their canonical projections 
$\bar{X}_1, \bar{X}_2, \ldots, \bar{X}_r \in \lie{g}_{r,c}^\CC/\gamma_2(\lie{g}_{r,c}^\CC)$ form a basis of $\lie{g}_{r,c}^\CC/\gamma_2(\lie{g}_{r,c}^\CC)$.
This implies that $X_1,X_2, \ldots, X_r$ are free generators of the free nilpotent Lie algebra $\lie{g}_{r,c}^\CC$. We can assume that $H$ is a Hall basis  with $H_1=\{X_1,X_2,\ldots,X_r\}$ and that $X_j$ is an eigenvector with eigenvalue $\lambda_j$.
By induction, it now follows that if $X\in H_i$ ($1\leq i \leq c$), then $X$ is an eigenvector for $d\varphi$ with eigenvalue $f_\lambda(X)$. Indeed, assume that $i\geq 2$ and the claim already holds for smaller values of $i$, then $X$ is of the form 
$X=[U,V]$ with $U\in H_k$ and $V\in H_l$ for some $k,l<i$. Then 
\[d\varphi(X)= d\varphi[U,V] =[d\varphi U, d\varphi V]= [f_\lambda(U) U, f_\lambda(V) V] =f_\lambda(U)f_\lambda(V)[U,V]=
f_\lambda(X) X.\]
\sloppy As the canonical projections of the vectors in $H_i$ form a basis for the vector space $\gamma_i(\lie{g}_{r,c}^\CC)/\gamma_{i+1}(\lie{g}_{r,c}^\CC)$, it follows that the collection of eigenvalues of $d\varphi_i$, and hence also of $\varphi_i$, is exactly the collection of 
values $f_\lambda(X)$, where $X$ ranges over all vectors in $H_i$.
\end{proof}

\begin{ex} Continuing \cref{exhall} and \cref{exlambda} we find that when $\lambda_1,\lambda_2, \ldots,
\lambda_r$ are the eigenvalues of $\varphi_1$, then the eigenvalues of $\varphi_2$ are 
\[ \lambda_i \lambda_j \mbox{ with } 1 \leq i < j \leq r\]
and those of $\varphi_3$ are
\[ \lambda_i \lambda_j \lambda_k\mbox{ with }1 \leq j < k \leq r \mbox{ and }1 \leq j \leq i \leq r.\]
\end{ex}

\medskip

As we will focus on diffeomorphisms, we are especially interested in the case that $\varphi$ is an automorphism (and not just any morphism). In this case the induced map $\varphi_1$ will be an automorphism of $\ZZ^r$.
 We can consider the morphism
\[ \psi:\Aut(N_{r,c}) \to \Aut(\ZZ^r): \varphi \mapsto \varphi_1\]
which is onto. Indeed, it is well known that the analogous map $\Aut(F_r) \to  \Aut(\ZZ^r)$ for the free group is onto (\cite[Theorem N4, Section 3.5]{mks76-1}). Since all automorphisms of $F_r$ induce an automorphism on $N_{r,c}$, the surjectivity of $\psi$ follows immediately.

\medskip

As explained above, $R(\varphi)$ only depends on the eigenvalues of $d\varphi$, which are completely determined by the eigenvalues of $\varphi_1$ (by \cref{lem:eigenvalues}). Hence, it is enough to know the characteristic polynomial of $\varphi_1$, which is of the form
\begin{equation}
\label{eq:degrpoly}
p(x) = x^r + a_{r-1}x^{r-1} + \cdots + a_1 x + a_0,
\end{equation}
where all $a_i\in\ZZ$ and $a_0=\pm 1$ (since $a_0=\pm \det(\varphi_1)$).

Conversely, any monic polynomial of degree $r$ of the form \eqref{eq:degrpoly} (still with $a_i\in \ZZ$ and $a_0=\pm 1$) is the characteristic polynomial of its companion matrix $C_p\in \GL_r(\ZZ)$, where   
\begin{equation*}
C_p = \begin{pmatrix}
	&			&	& -a_0\\
1	& 			&	& -a_1\\
	& \ddots	& 	& \vdots\\
	& 			& 1 & -a_{r-1}
\end{pmatrix}.
\end{equation*}

As $\psi$ is surjective, we know that there exists an automorphism $\varphi\in \Aut(N_{r,c})$ with $\varphi_1=C_p$.
So instead of focusing on the automorphisms $\varphi$, we will in the sequel focus on the corresponding characteristic polynomial. Let $p(x)$ be a polynomial of the form \eqref{eq:degrpoly}. We will denote by $R_c(p(x))$ the Reidemeister number of any automorphism $\varphi$ of $N_{r,c}$, such that the corresponding automorphism $\varphi_1$ has $p(x)$ as its characteristic polynomial. 

Thus, in order to calculate the Reidemeister spectrum of \(N_{r,c}\), we have to compute all possible numbers \(R_c(p(x))\) for all possible polynomials $p(x)$.

\section{The Reidemeister spectrum of \(N_{r,2}\)}

In this section we want to prove that the Reidemeister spectrum of $N_{r,2}$ is full for all $r\geq 4$ (and we will also compute the spectrum for $r=2$ and $r=3$). For this we need to study the numbers $R_2(p(x))$. 
So we consider any polynomial \(p(x)\) with integer coefficients of the form \eqref{eq:degrpoly} with \(a_0 = \pm 1\) and denote its roots (which are complex numbers in general) by \(\lambda_1, \lambda_2, \dots, \lambda_r\).  Then the Reidemeister number associated to this polynomial is 
\begin{equation}\label{eq:computeR2}
R_2(p(x)) = \left| \prod_{i=1}^r (1-\lambda_i) \prod_{1\leq i < j \leq r} (1-\lambda_i\lambda_j) \right|.
\end{equation}
For the sake of simplicity, we will always assume that \(p(x)\) is a polynomial such that \(R(p(x)) \neq \infty\) (i.e.\ none of $\lambda_i$ nor the $\lambda_i\lambda_j$ in the expression above equals 1). Note that the first factor of \(R_2(p(x))\) is exactly \(p(1) = \displaystyle\sum_{i=0}^{r-1}a_i + 1\).
\subsection{The case \(r=2\)}
We only have two roots \(\lambda_1, \lambda_2\), for which \(\lambda_1\lambda_2 = a_0\), hence
\begin{equation*}
R_2(p(x)) = \left| (a_0+a_1+1)(1-a_0)\right|.
\end{equation*}
Since we assume \(R_2(p(x)) \neq \infty\) we must have \(a_0 = -1\). Then \(R_2(p(x)) = 2|a_1|\), so only even numbers can be Reidemeister numbers. If \(q_n(x) = x^2 + nx -1\) with \(n \in \mathbb{N}\), then \(R_2(q_n(x)) = 2n\), hence \(\Spec_R(N_{2,2}) = 2\mathbb{N} \cup \{\infty\}\). The result of the computation coincides with \cite[Section 3]{roma11-1}, where  \(\Spec_R(N_{2,2})\) was computed via other techniques.

\subsection{The case \(r=3\)} There are three roots \(\lambda_1,\lambda_2,\lambda_3\) such that for \(i,j,k\) all distinct, \(\lambda_j\lambda_k = -a_0/\lambda_i\) and \(\lambda_i\lambda_j\lambda_k = a_0 = \pm 1\). Hence we may rewrite the Reidemeister number as 
\begin{align*}
R_2(p(x)) &= \left|\prod_{i=1}^3 (1-\lambda_i)\prod_{i=1}^{3}\left(1+\frac{a_0}{\lambda_i}\right)\right|\\
&= \left|\prod_{i=1}^3 (1-\lambda_i)\prod_{i=1}^3 (\lambda_i+a_0)\right|\\
&= |p(1)p(-a_0)|\\
&= \begin{cases}
(a_2+a_1)^2 & \textrm{ if } a_0 = -1\\
|(a_2+1)^2-(a_1+1)^2| & \textrm{ if } a_0 = 1.
\end{cases}
\end{align*}
Thus \(R_2(p(x))\) is a square or the difference of two squares, so it must be a multiple of four or an odd number. If \(q_n(x) = x^3 + nx^2 + (n-1)x +1\) with \(n \in \mathbb{N} \cup \{0\}\), then \(R_2(q_n(x)) = 2n+1\); and if \(r_n(x) = x^3 + n x^2 + (n-2)x +1\) with \(n \in \mathbb{N}\), then \(R_2(r_n(x)) = 4n\). Hence \(\Spec_R(N_{3,2}) = (2\mathbb{N}-1) \cup 4\mathbb{N} \cup \{\infty\}\). Again, this result coincides with that of  \cite[Section 3]{roma11-1}.
\subsection{Polynomials of even degree \(r \geq 4\)}
Let \(r = 2m\) and let \(n \in \NN\) be arbitrary. For the polynomial
\[ p_{2m,n}(x) = x^{2m} - x^{m+1} + (n-1) x^m + 1 \]
with roots \(\lambda_1,\lambda_2,\ldots,\lambda_{2m}\), we will show that \(R_2(p_{2m,n}(x)) = n\). 
This polynomial was first considered in \cite{mijl14-1} where it was also conjectured that indeed \(R_2(p_{2m,n}(x)) = n\). In her thesis  M.~Mijle checked this conjecture for $m=2,3, \ldots 9$.

In the computations below, we will simply write $p(x)$ instead of $p_{2m,n}(x)$. 
The first factor in the computation of $R_2(p(x))$, see \eqref{eq:computeR2}, is \(p(1) = n\), so it suffices to prove that
\begin{equation}
\label{eq:prodineqjis1}
\left[ \prod_{1 \leq i <  j \leq 2m} (1-\lambda_i \lambda_j) \right]^2 = \prod_{i\neq j} (1-\lambda_i \lambda_j) = 1. 
\end{equation}
We note that \(\prod_i \lambda_i = 1\) because \(p(x)\) has even degree and has constant term equal to \(1\). Also, if \(\lambda_i\) is a root of \(p(x)\) then
\[ \lambda_i^{2m} + (n-1) \lambda_i^m + 1 = \lambda_i^{m+1}, \]
so
\begin{align*}
\lambda_i^{2m} p\left(\frac{1}{\lambda_i}\right) &= \lambda_i^{2m} + (n-1)\lambda_i^m - \lambda_i^{m-1} + 1 \\
&= \lambda_i^{m+1} - \lambda_i^{m-1} \\
&= -\lambda_i^{m-1} (1-\lambda_i^2),
\end{align*}
giving
\begin{equation}
\label{eq:p1overlambda}
p\left(\frac{1}{\lambda_i}\right) = -\lambda_i^{-m-1}(1-\lambda_i^2).
\end{equation}
We want a polynomial whose roots include \(\lambda_i \lambda_j\), so to this end we define
\[ q(x) = \prod_{i = 1}^{2m} p\left(\frac{x}{\lambda_i}\right). \]
We then calculate
\begin{align}
q(x) &= \prod_{i=1}^{2m} p \left(\frac{x}{\lambda_i}\right) \nonumber\\
&= \prod_{i=1}^{2m} \prod_{j=1}^{2m} \left(\frac{x}{\lambda_i} - \lambda_j\right) \nonumber\\
&= \prod_{i=1}^{2m}  \frac{1}{\lambda_i^{2m}} \prod_{j=1}^{2m} (x-\lambda_i \lambda_j) \nonumber\\
&= \left[ \prod_{i=1}^{2m} \frac{1}{\lambda_i} \right]^{2m} \prod_{1 \leq i,j \leq 2m} (x-\lambda_i \lambda_j) \nonumber\\
&= \prod_{1 \leq i,j \leq 2m} (x-\lambda_i \lambda_j) \nonumber\\
&= \prod_{i\neq j} (x-\lambda_i \lambda_j)  \prod_{i=1}^{2m} (x-\lambda_i^2), \label{eq:qx}
\end{align}
where in the second-to-last line we used the noted fact that \(\prod_i \lambda_i = 1\). Consider
\begin{equation}
\label{eq:q1one}
q(1) = \prod_{i\neq j} (1-\lambda_i \lambda_j) \prod_{i=1}^{2m} (1-\lambda_i^2).
\end{equation}
For comparison, from \eqref{eq:p1overlambda} we obtain an alternate representation for \(q(1)\):
\begin{align}
q(1) &= \prod_{i=1}^{2m} p\left(\frac{1}{\lambda_i}\right) \nonumber\\
&= \prod_{i=1}^{2m} \left[-\lambda_i^{-m-1}(1-\lambda_i^2)\right] \nonumber\\
&= \left[\prod_{i=1}^{2m} \lambda_i \right]^{-m-1} \prod_{i=1}^{2m} (1-\lambda_i^2) \nonumber\\
&= \prod_{i=1}^{2m} (1-\lambda_i^2), \label{eq:q1two}
\end{align}
where we have again used the fact that \(\prod_i \lambda_i = 1\). Since \(n \in \NN\) and \(p(1) = n\), the number \(1\) is not a root of \(p\). 

If \(-1\) is not a root of \(p\), then \(\lambda_i^2 \neq 1\) for all \(i\), so the factor \(\prod_i (1-\lambda_i^2)\) in both \eqref{eq:q1one} and \eqref{eq:q1two} is non-zero. The desired identity \eqref{eq:prodineqjis1} then follows.

Now suppose that \(-1\) is a root of \(p\), then \(n = 2(-1)^{m+1}\) and 
\begin{align}
p'(-1) &=  2m(-1)^{2m-1} - (m+1)(-1)^m + m(n-1)(-1)^{m-1}  \nonumber\\
&= -2m + (m+1)(-1)^{m+1} + m(2(-1)^{m+1}-1)(-1)^{m-1}.\nonumber\\
&= (-1)^{m+1}. \label{eq:p'(-1)}
\end{align}
Hence \(-1\) is not a double root, so we can call this root \(\lambda_1\). For \(x \neq 1\) we can divide both sides of \eqref{eq:qx} by \(x-1\) to get
\[
\frac{q(x)}{x-1} = \prod_{i\neq j} (x-\lambda_i \lambda_j) \prod_{i=2}^{2m} (x-\lambda_i^2),
\]
and hence
\begin{equation}
\label{eq:qxoverx-1}
\lim_{x \to 1} \frac{q(x)}{x-1} = \prod_{i\neq j} (1-\lambda_i \lambda_j)  \prod_{i=2}^{2m} (1-\lambda_i^2).
\end{equation}
Alternatively,
\[
\frac{q(x)}{x-1} = \frac{p(-x)}{x-1} \prod_{i=2}^{2m} p\!\left(\frac{x}{\lambda_i}\right),
\]
so that
\begin{align}
\lim_{x \to 1} \frac{q(x)}{x-1} &= \left. \frac{d}{dx}p(-x) \right|_{x=1} \prod_{i=2}^{2m} p\!\left(\frac{1}{\lambda_i}\right) \nonumber\\
&= - p'(-1) \prod_{i=2}^{2m} \left[-\lambda_i^{-m-1}(1-\lambda_i^2)\right] \nonumber\\
&=  p'(-1) \left(\prod_{i=2}^{2m} \lambda_i \right)^{-m-1} \prod_{i=2}^{2m} (1-\lambda_i^2) \nonumber\\
&= \prod_{i=2}^{2m} (1-\lambda_i^2) \label{eq:qxoverx-1two}
\end{align}
where in the second line we used identity \eqref{eq:p1overlambda} and in the last line we used both \eqref{eq:p'(-1)} and the fact that
\[ \prod_{i=2}^{2m} \lambda_i = \lambda_1^{-1} = \lambda_1 = -1. \]
By comparing \eqref{eq:qxoverx-1} and \eqref{eq:qxoverx-1two} the desired identity \eqref{eq:prodineqjis1} follows. 

\medskip

As a conclusion of this computation we find:

\begin{prop}\label{prop:full1}
Let $m\geq 2$ be an integer, then  \(\Spec_R(N_{2m,2})\) is full. 
\end{prop}

\subsection{Polynomials of odd degree \(r \geq 5\)} Let \(r = 2m+1\) and let \(n \in \NN\) be arbitrary. For the polynomial
\[ p_{2m+1,n}(x) = x^{2m+1} + (n+1)x^{m+2} + (1-n)x^{m+1} + (n-1)x^m - nx^{m-1}-1, \]
with roots \(\lambda_1,\lambda_2,\ldots,\lambda_{2m+1}\), we will show that \(R_2(p_{2m+1,n}(x)) = n+c(m)\) with
\begin{equation*}
c(m) = 2+\cos\left(m\frac{\pi}{3}\right)+\sqrt{3}\sin\left(m\frac{\pi}{3}\right) = \begin{cases}
0 & \text{ if } m \equiv 4 \pmod 6\\
1 & \text{ if } m \equiv 3 \pmod 6 \text{ or } m \equiv 5 \pmod 6\\
3 & \text{ if } m \equiv 0 \pmod 6 \text{ or } m \equiv 2 \pmod 6\\
4 & \text{ if } m \equiv 1 \pmod 6\\
\end{cases}.
\end{equation*} 
It then follows that \(R_2(p_{2m+1,n-c(m)}(x)) = n\). The proof uses similar techniques as for the case where \(r\) is even. Again, during the computations, we will simply write $p(x)$ instead of $p_{2m+1,n}(x)$.

Now, the first factor in \eqref{eq:computeR2} is \(p(1) = 1\), so it suffices to prove that
\begin{equation}
\label{eq:prodineqjis1odd}
\left|\prod_{1 \leq i <  j \leq 2m+1} (1-\lambda_i \lambda_j)\right| =  n+c(m).
\end{equation}
We first calculate some specific values of \(p(x)\):
\begin{align}
p(1) &= \prod_{i=1}^{2m+1}(1-\lambda_i)&&= 1,\label{eq:p(1)two}\\
p(-1) &=-\prod_{i=1}^{2m+1}(1+\lambda_i)&&= (-1)^m(4n-1)-2.\label{eq:p(-1)two}
\end{align}
We find that \(1\) and \(-1\) both are not roots of \(p(x)\). Also, for any root \(\lambda_i\) we have
\begin{equation*}
-\lambda_i^{2m+1} = (n+1)\lambda_i^{m+2} + (1-n)\lambda_i^{m+1} + (n-1)\lambda_i^m - n\lambda_i^{m-1} -1,
\end{equation*}
so
\begin{align*}
\lambda_i^{2m+1}p\left(\frac{1}{\lambda_i}\right) &= -\lambda_i^{2m+1} -n\lambda_i^{m+2} + (n-1)\lambda_i^{m+1} +(1-n)\lambda_i^m + (n+1)\lambda_i^{m-1} +1\\
&= \lambda_i^{m+2} + \lambda_i^{m-1}\\
&= \lambda_i^{m-1}\left(  1 + \lambda_i^{3}\right)\\
&= \lambda_i^{m-1}\left(  1 + \lambda_i\right)\left(  e^{\frac{\pi}{3}i} - \lambda_i\right)\left(  e^{-\frac{\pi}{3}i} - \lambda_i\right),
\end{align*}
giving
\begin{equation}
\label{eq:p1overlambdatwo}
p\left(\frac{1}{\lambda_i}\right) = \lambda_i^{-m-2}\left(  1 + \lambda_i\right)\left(  e^{\frac{\pi}{3}i} - \lambda_i\right)\left(  e^{-\frac{\pi}{3}i} - \lambda_i\right).
\end{equation}
Again, we define a new polynomial \(q(x)\) as
\begin{equation*}
q(x) = \prod_{i = 1}^{2m+1}p\left(\frac{x}{\lambda_i}\right),
\end{equation*}
and once again
\begin{equation*}
q(x) = \prod_{i\neq j} (x-\lambda_i \lambda_j)  \prod_{i=1}^{2m+1} (x-\lambda_i^2).
\end{equation*}
Let us evaluate \(q(x)\) at \(x = 1\):
\begin{align}
q(1) &= \prod_{i \neq j}\left(1-\lambda_i\lambda_j\right)\prod_{i  = 1}^{2m+1}\left(1-\lambda_i^2\right)\nonumber\\
&= \prod_{i \neq j}\left(1-\lambda_i\lambda_j\right)\prod_{i  = 1}^{2m+1} \left(1-\lambda_i\right)\prod_{i  = 1}^{2m+1}\left(1+\lambda_i\right)\nonumber\\
&= -p(-1)\prod_{i \neq j}\left(1-\lambda_i\lambda_j\right),\label{eq:q1oneodd}
\end{align}
where we used \eqref{eq:p(1)two} and \eqref{eq:p(-1)two} in the last step.

We evaluate \(q(x)\) at \(x=1\) using \eqref{eq:p1overlambdatwo}:
\begin{align}
q(1) &= \prod_{i=1}^{2m+1}p\left(\frac{1}{\lambda_i}\right)\nonumber\\ 
&= \prod_{i = 1}^{2m+1}\lambda_i^{-m-2}\left(  1 + \lambda_i\right)\left(  e^{\frac{\pi}{3}i} - \lambda_i\right)\left(  e^{-\frac{\pi}{3}i} - \lambda_i\right)\nonumber\\
&= \left[ \prod_{i = 1}^{2m+1}\lambda_j\right]^{-m-2} \prod_{i = 1}^{2m+1} \left(  1 + \lambda_i\right)\prod_{i = 1}^{2m+1}\left(  e^{\frac{\pi}{3}i} - \lambda_i\right)\prod_{i = 1}^{2m+1}\left(  e^{-\frac{\pi}{3}i} - \lambda_i\right)\nonumber\\
&= -p(-1) p\left(e^{\frac{\pi}{3}i}\right)p\left(e^{-\frac{\pi}{3}i}\right)\nonumber\\
&= -p(-1) \left|p\left(e^{\frac{\pi}{3}i}\right)\right|^2,\label{eq:q1twoodd}
\end{align}
where we used that \(p\left(\overline{z}\right) = \overline{p\left(z\right)}\) for any \(z \in \mathbb{C}\), since \(q(x)\) only has real coefficients. Comparing equations \eqref{eq:q1oneodd} and \eqref{eq:q1twoodd} now gives 
\begin{equation*}
\prod_{i \neq j}\left(1-\lambda_i\lambda_j\right) = \left|p\left(e^{\frac{\pi}{3}i}\right)\right|^2,
\end{equation*}
or by using that the product is symmetric in the indices, that
\begin{equation*}
\left|\prod_{1 \leq i < j \leq 2m+1}\left(1-\lambda_i\lambda_j\right)\right| = \left|p\left(e^{\frac{\pi}{3}i}\right)\right|.
\end{equation*}
One may now evaluate \(\left|p\left(e^{\frac{\pi}{3}i}\right)\right|\) to obtain
\begin{equation*}
\left|p\left(e^{\frac{\pi}{3}i}\right)\right| = \left|n+2+\cos\left(m\frac{\pi}{3}\right)+\sqrt{3}\sin\left(m\frac{\pi}{3}\right)\right| = n+c(m).
\end{equation*}
This proves the following proposition

\begin{prop}\label{prop:full2}
Let $m\geq 2$ be an integer, then  \(\Spec_R(N_{2m+1,2})\) is full. 
\end{prop}

\medskip

As a conclusion we have:
\begin{thm} Let $M$ be a nilmanifold whose fundamental group is the free 2-step nilpotent group $N_{r,2}$ with $r\geq 4$. Then, for any non-negative integer $n$ there exists a self diffeomorphism $h_n$ of $M$ such that $h_n$ has exactly $n$ fixed points and any self-map $f$ of $M$ which is homotopic to $h_n$ has at least $n$ fixed points.
\end{thm}
\begin{proof}
By \cref{prop:full1} and \cref{prop:full2} we know that  \(\Spec_R(N_{r,2})\) is full for $r\geq 4$. The theorem then follows from 
\cref{cor:main}.
\end{proof}

\section{The Reidemeister spectrum of \(N_{r,c}\) with \(c > 2\)}
A symmetric polynomial with coefficients in \(\ZZ\) and variables \(\lambda_1, \dots, \lambda_r\), can be expressed in terms of the elementary symmetric polynomials. If the \(\lambda_i\) are then interpreted as the roots of a polynomial \(p(x) = \sum_{i=0}^r a_ix^i\), then these elementary symmetric polynomials are, up to sign, the coefficients \(a_i\). In other words:
\begin{equation*}
q \in \ZZ[\lambda_1, \dots, \lambda_r] \text{ and } q \text{ symmetric } \implies q \in \mathbb{Z}[a_0, \dots, a_r].
\end{equation*}
To calculate the Reidemeister spectrum of \(N_{r,c}\) with \(c > 2\), we will adopt a ``divide and conquer'' strategy, splitting up  \(R_c(p(x))\) in factors that are each symmetric polynomials in the roots \(\lambda_i\), and calculating these factors in terms of the coefficients \(a_i\).

Let us calculate the spectrum of \(N_{2,3}\) and \(N_{3,3}\) to demonstrate this approach. Let \(p(x)\) be of the form \eqref{eq:degrpoly} again with \(a_0 = \pm 1\), then
\begin{equation*}
R_3(p(x)) = \left| \prod_{i=1}^r (1-\lambda_i) \prod_{i < j} (1-\lambda_i\lambda_j) \prod_{\substack{j<k \\ j \leq i}}(1-\lambda_i\lambda_j\lambda_k) \right|.
\end{equation*}
Since we assume that \(R_3(p(x)) \neq \infty\) and $r=2$ or $3$, we know that \(\prod \lambda_i = (-1)^r a_0 = -1\).

\medskip

 For \(N_{2,3}\), we already know the first two factors from the calculations we made for \(N_{2,2}\), they are \(a_1\) and \(2\) respectively. The third factor becomes
\begin{equation*}
\prod_{\substack{j< k \\ j \leq i }}(1-\lambda_i\lambda_j\lambda_k) = (1-\lambda_1^2\lambda_2)(1-\lambda_1\lambda_2^2) = (1+\lambda_1)(1+\lambda_2) = p(-1) = a_1,
\end{equation*}
so all factors combined give
\begin{equation*}
R_3(p(x)) = 2a_1^2.
\end{equation*}
So \(R_3(p(x) )\) must be two times a square, and for the polynomials \(q_n = x^2 + nx -1\) with \(n \in \NN\) we have \(R_3(q_n(x)) = 2n^2\), hence \(\Spec_R(N_{2,3}) = 2\NN^2 \cup \{\infty\}\). This result was also obtained, using another approach,  in 
\cite[Section 3]{roma11-1}.

\medskip

For \(N_{3,3}\), we already know the first two factors from the calculations we did for \(N_{3,2}\):
\begin{equation*}
\prod_{i=1}^3(1-\lambda_i) = p(1) = 2 + a_1 + a_2, \quad \prod_{i<j}(1-\lambda_i\lambda_j) = -p(-1) = a_1 - a_2.
\end{equation*}
We can split the third factor in two symmetric polynomials:
\begin{equation*}
\prod_{\substack{j< k \\ j\leq i }}(1-\lambda_i\lambda_j\lambda_k) = (1-\lambda_1\lambda_2\lambda_3)^2 \prod_{i \neq j}(1-\lambda_i^2\lambda_j).
\end{equation*}
The first factor is clearly \(2^2 = 4\), and using that \(\lambda_1\lambda_2\lambda_3 = -1\) and \(\lambda_1+\lambda_2+\lambda_3 = -a_2\), we obtain for the second factor
\begin{align*}
\prod_{i \neq j}(1-\lambda_i^2\lambda_j) &= \prod_{i \neq j \neq k \neq i} \left(1-\frac{\lambda_i^2\lambda_j\lambda_k}{\lambda_k}\right)\\ &= \prod_{i \neq k}\left(1+\frac{\lambda_i}{\lambda_k}\right)\\
&=  \left(\prod_{k=1}^3\lambda_k\right)^{-2} \prod_{i \neq k}(\lambda_i + \lambda_k)\\
&=  \left[\prod_{j=1}^3(-a_2-\lambda_j)\right]^2\\
&= p(-a_2)^2\\
&= (1-a_1a_2)^2,
\end{align*}
so by putting everything together we find 
\begin{equation*}
R_3(p(x)) = 4\left|(2+a_1+a_2)(a_1-a_2)(1-a_1a_2)^2\right|.
\end{equation*}
Substituting \(a = 1+a_1\) and \(b = 1+a_2\), we may rewrite this as
\begin{equation*}
R_3(p(x)) = 4|a^2 - b^2|(a + b - ab)^2,
\end{equation*}
and in particular, if \(a \neq \pm b\), then  \(|a^2 - b^2| \geq |a| + |b|\), hence \(R_3(p(x)) \geq 4(|a|+|b|)\). This allows us to, in some sense, calculate the Reidemeister spectrum of \(N_{3,3}\). By calculating \(R_3(p(x))\) for all couples \((a,b)\) with \(|a|,|b| \leq 250\), we know that the Reidemeister numbers less than \(1000\) are exactly \(4\), \(12\), \(20\), \(32\), \(60\), \(64\), \(96\), \(108\), \(140\), \(192\), \(252\), \(300\), \(320\), \(324\), \(396\), \(480\), \(500\), \(572\), \(672\), \(700\),
\(756\), \(780\), \(800\), \(896\) and  \(980\).

To give a general idea on what numbers can be in the spectrum, consider the different values of \(a\) and \(b \mod 2\):
\begin{itemize}
	\item \(a,b \equiv 0 \mod 2\). Then \(|a^2 - b^2|\) is a multiple of \(4\) and \(a+b-ab\) is a multiple of \(2\). Hence \(R_3(p(x)) \in 64\NN\). 
	\item \(a \equiv 0, b \equiv 1 \mod 2\) or vice versa. Then both \(|a^2 - b^2|\) and \(a + b - ab\) are odd, hence \(R_3(p(x)) \in 4(2\NN-1)\). 
	\item \(a,b \equiv 1 \mod 2\). Then \(|a^2 - b^2|\) is a multiple of \(8\) and \(a+b-ab\) is odd. Hence \(R_3(p(x)) \in 32\NN\). 
\end{itemize}
Together with the calculated Reidemeister numbers mentioned earlier, we may then state that \(\Spec_R(N_{3,3}) \subsetneq 32\NN \cup 4(2\NN-1) \cup \{\infty\}\).

\medskip

A similar approach can be done for all \(N_{r,c}\), though for \(r > 3\) the calculations quickly become rather bothersome.

\medskip

Of particular interest is the following result:
\begin{prop}
	Let \(N_{r,c}\) be a free nilpotent group with \(c \geq r\). Then the Reidemeister spectrum of \(N_{r,c}\) is not full.
\end{prop}
\begin{proof}
	For any polynomial \(p(x)\), the \(r\)-th factor of \(R_c(p(x))\) will be a product of the form
	\begin{equation*}
	\prod (1-\lambda_{i_1}\lambda_{i_2} \cdots \lambda_{i_r}).
	\end{equation*}
	Splitting this further up in factors that are each symmetric polynomials, one of them will be \begin{equation*}
	(1-\lambda_1\lambda_2\cdots\lambda_r) = 1-a_0,
	\end{equation*}
	which is \(0\) or \(2\) depending on the constant term \(a_0\). Hence either \(R_c(p(x)) = \infty\) or  else \(2\) divides \(R_c(p(x))\), and therefore \(\Spec_R(N_{r,c}) \subseteq 2\NN \cup \{\infty\}\). 
\end{proof}
Based on some experimenting with polynomials with small coefficients, we suspect the following conjecture is true:
\begin{con}
	Let \(N_{r,c}\) be a free nilpotent group with \(r \geq 2\) and  \(c>2\). Then the Reidemeister spectrum of \(N_{r,c}\) is not full.
\end{con}

In this context it is also useful to recall the main result of \cite[Theorem 2.5]{dg14-1}:
\begin{thm}
Let $r\geq 2$, then $N_{r,c}$ has the $R_\infty$ property if and only if $c \geq 2r$.
\end{thm}

\medskip

\section{The extended Reidemeister spectrum of torsion-free nilpotent groups}
Similarly to the Reidemeister spectrum, we can define the \emph{extended Reidemeister spectrum} of a group \(G\) as
\begin{equation*}
\ESpec_R(G) = \{ R(\varphi) \mid \varphi \in \End(G)\}.
\end{equation*}
Let \(N\) be any torsion-free, finitely generated, nilpotent group (the free nilpotent groups are examples of this). Then \(N\) is a poly-\(\ZZ\) group and hence \(N \cong M \rtimes \ZZ\) for some normal  subgroup \(M \lhd N\). Consider the endomorphism 
\begin{equation*}
\varphi_n: M \rtimes \ZZ \to M \rtimes \ZZ: (m,z) \mapsto (1,nz).
\end{equation*}
It is easy to see that for $n\geq 2$ \(R(\varphi_n) = n-1\), hence \(\ESpec_R(N)\) is full.

As a consequence of the fact that \(\ESpec_R(N)\) is full for any finitely generated torsion free nilpotent group and using exactly the same proofs as for \cref{thm:main} and \cref{cor:main}, we now find the following:

\begin{thm}
Let $M$ be any nilmanifold. For any nonnegative integer $n$ there exists a smooth self-map $h_n$ of $M$ with exactly $n$ fixed points. Moreover, any selfmap $f$ of $M$ which is homotopic to $h_n$ has at least $n$ fixed points.
\end{thm}

\printbibliography

\end{document}